\documentclass[graybox]{svmult}

\usepackage{mathptmx}        
\usepackage{helvet}          
\usepackage{courier}         

\usepackage{makeidx}         
\usepackage{graphicx}        
\usepackage{multicol}        
\usepackage[bottom]{footmisc}

\usepackage{url}             

\makeindex             

\begin{document}




\mainmatter
%
%
%
%

{ 

\title*{%
    Time integration as filtering:
    a space--time discretization-aware
    LES formulation
}
\titlerunning{Space--time discretization-aware LES}

\author{Syver D{\o}ving Agdestein}

\institute{%
    Centrum Wiskunde \& Informatica,
    Science Park 123,
    1098 XG Amsterdam,
    The Netherlands,
    \email{sda@cwi.nl}
}

\maketitle

\abstract{%
Discretization-aware LES yields an \emph{exact} expression for the discrete
target flux in finite-volume LES by recognizing that a coarse finite difference
is a top-hat-filtered exact derivative (the ``filter-swap'' property).
That argument is purely spatial; here we observe that the forward-Euler
time difference is itself a (one-sided) top-hat-filtered exact \emph{time}
derivative, and repeat the construction in space--time. The resulting
exact discrete flux decomposition extends the spatial one with a single
temporal term: a flux-quadrature error that shrinks with the quadrature
order of the time integrator. In a Burgers experiment this term grows with
the CFL number while the spatial terms do not, and a Smagorinsky closure
augmented with its leading order---a Lax--Wendroff-type diffusion---stays
accurate at coarse time steps where space-only closures degrade.
}

\section{Introduction}
\label{ag:sec:intro}
The turbulence closure problem is usually posed for the \emph{spatially}
filtered equations. In the implicit-LES setting the only spatial filter is the
grid (top-hat) filter itself, of width $h > 0$, and the residual is dominated
by discretization effects---the
numerical flux and the discrete divergence---rather than by a sub-filter stress
alone \cite{chow-moin}. Discretization-aware LES
\cite{agdestein-jcp} addresses this by deriving an \emph{exact} residual
stress, hinging on the
\emph{filter-swap} property: the coarse central finite difference
$\partial^h_x u(x) \mathrel{\mathop:}=
(u(x + h/2) - u(x - h/2)) / h$
equals a top-hat-filtered exact derivative,
\begin{equation}
\label{ag:eq:space-swap}
\partial^h_x = \overline{(\cdot)}^{h} \circ \partial_x ,
\end{equation}
where $\overline{(\cdot)}^{h}$ is the volume-averaging (top-hat) grid filter of width $h$
\cite{schumann,agdestein-jcp}. This argument is entirely spatial: the time
derivative $\partial_t$ is kept
exact, and the time integrator is applied only afterwards; the
time-discretization error therefore lives \emph{outside} the exact
residual-stress guarantee.

Filtering in \emph{time} is not in itself new: temporal LES models the
residual stress of an explicitly chosen, smooth causal time filter
\cite{pruett}, and the view of the LES equations as the result of two
operators---a smooth filter \emph{and} a discretization operator---goes back to \cite{carati-winckelmans}.
Our observation is different: the forward-Euler time integrator \emph{is} the
(one-sided top-hat) time filter, so the filter-swap is an \emph{exact} identity
and the only new contribution is a single quadrature term---no separate
temporal closure or commutation term arises.

\section{From spatial to space--time filter-swap}
\label{ag:sec:theory}
We consider the scalar conservation law
\begin{equation} \label{ag:eq:cl}
    \partial_t u + \partial_x f(u) = 0,
\end{equation}
with flux $f$; the experiments use the viscous Burgers flux
$f(u) \mathrel{\mathop:}= u^2 / 2 - \nu\, \partial_x u$.

\runinhead{Spatial filter-swap.}
Let $f^h$ be a grid-compatible numerical flux. Filtering
(\ref{ag:eq:cl}) and absorbing the grid filter into the spatial derivative
through (\ref{ag:eq:space-swap}) gives the spatially discrete, time-continuous
law of \cite{agdestein-jcp},
\begin{equation} \label{ag:eq:space-rst}
    \partial_t \bar{u}^{h} + \partial^h_x f(u) = 0 .
\end{equation}
The grid difference $\partial^h_x$ acts on the \emph{unfiltered} flux $f(u)$:
this is the content of the filter-swap, and it is exact. Replacing the
inaccessible $f(u)$ by the numerical flux $f^h(\bar{u}^h)$ recovers the usual
semi-discrete finite-volume scheme; the time derivative is still exact.

\runinhead{Temporal filter-swap.}
Let $\tau$ be the time step. The forward top-hat time filter and the forward
time difference,
\begin{equation} \label{ag:eq:time-swap}
    \overline{u}^{\tau}(t) \mathrel{\mathop:}= \frac{1}{\tau}\int_0^{\tau} u(t+s)\,\mathrm{d}s,
    \qquad
    \partial^\tau_t u(t) \mathrel{\mathop:}= \frac{u(t+\tau)-u(t)}{\tau}
    = \overline{\partial_t u}^{\,\tau}(t),
\end{equation}
satisfy the temporal filter-swap
$\partial^\tau_t = \overline{(\cdot)}^{\tau} \circ \partial_t$ exactly---the
temporal analogue of (\ref{ag:eq:space-swap}).
Forward Euler is this one-sided, first-order filter.

\runinhead{Space--time filter-swap.}
Applying the space--time filter
$\overline{(\cdot)}^{h, \tau} \mathrel{\mathop:}= \overline{(\cdot)}^{\tau} \circ \overline{(\cdot)}^{h}$
to (\ref{ag:eq:cl}) and absorbing each filter into its own derivative through
(\ref{ag:eq:space-swap}) and (\ref{ag:eq:time-swap}) yields a law that is
discrete in space \emph{and} time,
\begin{equation} \label{ag:eq:st-filtered}
    \partial^\tau_t\, \bar{u}^{h} + \partial^h_x\, \overline{f(u)}^{\tau} = 0 .
\end{equation}
This law is
exact: injecting $\overline{f(u)}^{\tau}$ as the flux makes the forward-Euler
scheme reproduce the grid-filtered DNS $\bar{u}^{h}$ at the time nodes with zero
residual. The only departure from the spatial law (\ref{ag:eq:space-rst}) is
that the flux is \emph{time-averaged over the step}; the resolved field stays
$\bar{u}^{h}$.
In practice $\overline{f(u)}^{\tau}$ is replaced by an instantaneous,
step-aware numerical flux
$f^{h, \tau}(\bar{u}^h) \approx \overline{f(u)}^{\tau}$ that uses only the
current field but may depend on $\tau$, giving the fully discrete scheme
$v^{n+1}_i \mathrel{\mathop:}= v^n_i - (\tau / h) [
f^{h, \tau}_{i + 1/2}(v^n) - f^{h, \tau}_{i - 1/2}(v^n) ]$
with $v^n_i \approx \bar{u}^{h}(x_i, t_n)$.
The $\tau$-awareness comes through a closure $m^{h, \tau}$ as $f^{h, \tau} \mathrel{\mathop:}= f^h + m^{h, \tau}$.

\section{Decomposition and modeling of the space--time residual flux}
\label{ag:sec:decomp}
The closure target is the \emph{residual flux}---the gap between the exact
space--time flux of (\ref{ag:eq:st-filtered}) and the numerical one,
\begin{equation} \label{ag:eq:residual}
    r^{h, \tau}(u) \mathrel{\mathop:}= \overline{f(u)}^{\tau} - f^h(\bar{u}^h) .
\end{equation}
Inserting the intermediate fluxes $f(u)$, $\overline{f(u)}^h$, and $f(\bar{u}^h)$
splits it into four errors,
\begin{eqnarray}
    r^{h, \tau} & = & r^h_{\mathrm{LES}} + r^h_{\mathrm{num}} + r^h_{\mathrm{div}} + r^\tau_{\mathrm{time}}, \label{ag:eq:decomp}\\
    r^h_{\mathrm{LES}}(u) & \mathrel{\mathop:}= & \overline{f(u)}^h - f(\bar{u}^h),
    \qquad
    r^h_{\mathrm{num}}(u) \mathrel{\mathop:}= f(\bar{u}^h) - f^h(\bar{u}^h), \nonumber\\
    r^h_{\mathrm{div}}(u) & \mathrel{\mathop:}= & f(u) - \overline{f(u)}^h,
    \qquad\;\;\,
    r^\tau_{\mathrm{time}}(u) \mathrel{\mathop:}= \overline{f(u)}^{\tau} - f(u) . \nonumber
\end{eqnarray}
The first three are the spatial terms of \cite{agdestein-jcp}:
$r^h_{\mathrm{LES}}$ is the classical LES commutator between grid filtering and
the nonlinear flux; $r^h_{\mathrm{num}}$ is the numerical-flux error; and
$r^h_{\mathrm{div}}$ corrects
$r^h_{\mathrm{LES}}$ for the \emph{discrete} divergence. The fourth term
$r^\tau_{\mathrm{time}}$ is new: a pure \emph{flux
time-quadrature} error---the step-averaged flux minus its instantaneous
value---that shrinks with the order of the time integrator (forward Euler is
the left-rectangle rule).

\runinhead{A space--time discretization-aware eddy-viscosity model.}
In \cite{agdestein-jcp} the three spatial terms are modeled together by a
discretization-aware Smagorinsky \cite{smagorinsky} term $m^h$ embedded in
the numerical flux, with
$f(u) \approx f^h(\bar{u}^h) + m^h(\bar{u}^h)$,
\begin{equation} \label{ag:eq:smag-h}
    m^h(\cdot) \mathrel{\mathop:}= -\, \theta^2 h^2\, |\partial^h_x \cdot|\, \partial^h_x \cdot ,
\end{equation}
whose single constant $\theta$ is tuned by dissipation matching against the
exact flux $f(u)$.
Expanding the
step average and using the conservation law (\ref{ag:eq:cl}) gives, for the
convective part of the flux,
\begin{equation} \label{ag:eq:lw}
    r^\tau_{\mathrm{time}}(u)
    = \frac{\tau}{2}\, \partial_t f(u) + O(\tau^2)
    = -\, \frac{\tau}{2}\, (f'(u))^2\, \partial_x u + O(\tau^2) ,
\end{equation}
a diffusive flux whose effective viscosity $\frac{\tau}{2} (f'(u))^2$ is
non-negative \emph{by construction}---precisely the Lax--Wendroff diffusion
that stabilizes forward-Euler convection \cite{lax-wendroff}. For the
space--time target $\overline{f(u)}^{\tau}$
we therefore augment the spatial closure with this term,
\begin{equation} \label{ag:eq:smag-ht}
    m^{h, \tau}(\cdot)
    \mathrel{\mathop:}= -\, \theta_h^2 h^2\, |\partial^h_x \cdot|\, \partial^h_x \cdot
    \, -\, \theta_\tau\, \tau\, (f'(\cdot))^2\, \partial^h_x \cdot ,
\end{equation}
with the theoretical value $\theta_\tau = \frac{1}{2}$ from
(\ref{ag:eq:lw}).

\section{Numerical experiment: space--time LES for Burgers}
\label{ag:sec:experiment}
We test the framework on viscous Burgers turbulence.
\footnote{The code is archived at \url{https://zenodo.org/records/20714279}.}
All results are
ensemble averages over $100$ random initial fields with the
smooth energy spectrum
$E(k) \propto (k / 10)^4 \mathrm{e}^{-2 (k / 10)^2}$, evolved on a periodic
domain of length $2\pi$ with
$\nu = 5 \times 10^{-4}$. The DNS is a second-order finite-volume scheme with
$13\,500$ cells, integrated with forward Euler
($\Delta t = 2.5 \times 10^{-5}$) to $t = 0.1$, when the solutions
have steepened into sawtooth profiles. Each coarse cell is the union of $45$
fine cells ($n = 300$; the odd ratio makes coarse and fine cell centers
coincide), and the coarse step is $\tau = M \Delta t$ with
$M \in \{25, 50, 100, 200, 400\}$, giving CFL numbers
$\tau \max_i |\bar{u}^h_i| / h$ between $0.06$ and $0.89$. With the coarse
faces a subset of the fine faces and $\tau$ a multiple of $\Delta t$,
both filter-swaps hold exactly in the fully discrete setting: forward Euler
driven by the recorded target flux $\overline{f(u)}^{\tau}$ reproduces
$\bar{u}^h$ to machine precision at every step size. Two coarse numerical fluxes
$f^h$ are considered: the energy-conserving central flux
\cite{agdestein-jcp} and its counterpart with first-order upwind convection.

\begin{figure}[t]
    \sidecaption
    \includegraphics[width=0.62\textwidth]{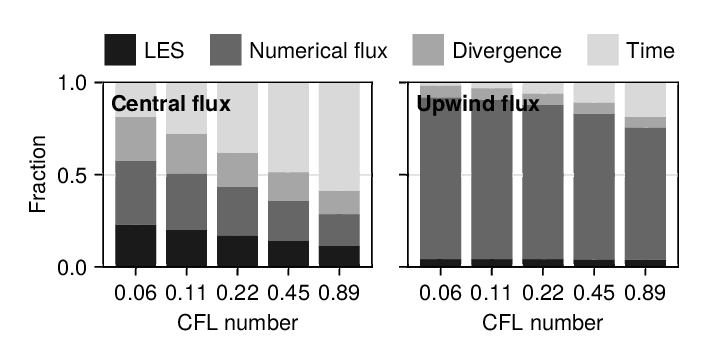}
    \caption{%
        Share $\langle | r_k | \rangle / \sum_l \langle | r_l | \rangle$
        of the four terms of the decomposition (\ref{ag:eq:decomp}) in the
        residual flux, versus CFL number, for the central (left) and upwind
        (right) numerical flux. The temporal share grows with the CFL number,
        while the spatial shares keep fixed proportions.
    }
    \label{ag:fig:fractions}
\end{figure}

\runinhead{A-priori analysis of the space--time flux.}
Figure~\ref{ag:fig:fractions} shows the share of each term of
(\ref{ag:eq:decomp}) in the residual flux, measured at the coarse faces at
$t = 0.1$. For the central flux the three spatial terms keep essentially
fixed proportions, while the temporal share grows from $0.19$ at CFL $0.06$
to $0.59$ at CFL $0.89$: at practical CFL numbers the flux time-quadrature
error is the \emph{largest single term}---precisely the term that
space-only formulations discard. For the upwind flux
the numerical-flux error dominates ($0.72$--$0.88$), but the temporal share
again grows steadily with the CFL number (from $0.02$ to $0.19$).

The corresponding dissipation rates
$D_k \mathrel{\mathop:}= -\langle r_k \, \partial^h_x \bar{u}^h \rangle$
(not shown) complement this picture. The LES and numerical-flux terms supply
nearly all of the dissipation, independently of
the CFL number, and the temporal
dissipation grows \emph{linearly} in $\tau$---the measured
$D_{\mathrm{time}} / \tau$ is constant to within a few percent over the
sixteenfold range of step sizes---exactly the Lax--Wendroff scaling of
(\ref{ag:eq:lw}), and the basis for fitting $\theta_\tau$. For the upwind
flux the numerical-flux term is strongly \emph{anti}-dissipative: the upwind
$f^h$ dissipates more than the exact flux, so the residual must return
energy. The temporal term, finally, while carrying nearly $60\,\%$
of the residual magnitude at CFL $0.89$, carries less than $10\,\%$ of the
dissipation: much of the time-quadrature error is non-dissipative and hence
invisible to an eddy-viscosity model.

\begin{figure}[t]
    \sidecaption
    \includegraphics[width=0.62\textwidth]{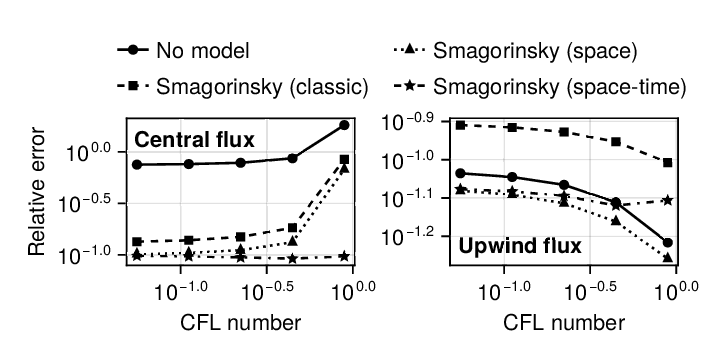}
    \caption{%
        Ensemble relative error of the coarse forward-Euler LES at the final time,
        versus CFL number, for the four closures, with the central (left) and upwind
        (right) numerical flux.
        Closures distinguished by line style and marker.
    }
    \label{ag:fig:error}
\end{figure}

\runinhead{A-posteriori analysis of the space--time-aware closure.}
We now run the coarse forward-Euler scheme of Sect.~\ref{ag:sec:theory} with
$f^{h, \tau} \mathrel{\mathop:}= f^h + m$ for four closures $m$:
\emph{no-model} ($m \mathrel{\mathop:}= 0$);
\emph{classical Smagorinsky} (\ref{ag:eq:smag-h}), $\theta$
fitted to $r^h_{\mathrm{LES}}$ alone; \emph{space-aware Smagorinsky}
(\ref{ag:eq:smag-h}), $\theta$ fitted to the three spatial terms; and
\emph{space--time-aware Smagorinsky} (\ref{ag:eq:smag-ht}),
$(\theta_h, \theta_\tau)$ fitted to the full $r^{h, \tau}$. All coefficients
are fitted a priori by dissipation matching ($\theta_\tau$ by least squares
over the five step sizes). This gives
$\theta^2_h = 0.202$ (central) and $-0.045$ (upwind); the fitted
$\theta_\tau = 0.21$ has the order of magnitude of
the theoretical $\frac{1}{2}$ of (\ref{ag:eq:lw}), reduced by the
coarse-graining $u \rightarrow \bar{u}^h$ in (\ref{ag:eq:smag-ht}).

Figure~\ref{ag:fig:error} shows the resulting ensemble relative error with
respect to $\bar{u}^h$ at $t = 0.1$. With
the central flux the no-model run is strongly under-dissipative (errors
above $0.75$) and becomes unstable at CFL $0.89$. The classic and
space-aware closures help at small steps ($0.13$ and $0.10$ at CFL $0.06$)
but, being step-agnostic, cannot supply the dissipation the time integrator
demands: their errors grow to $0.85$ and $0.68$ at CFL $0.89$. The
space--time-aware closure does supply it, and its error is essentially
\emph{independent of the CFL number}---between $0.092$ and $0.098$ across
the entire range: the single temporal term of
the decomposition is exactly what the space-only closures were missing.

The upwind flux reverses the picture. The numerical flux alone already
over-dissipates, so classical Smagorinsky---blind to the discretization---is
\emph{worse than no model} at every step size, while the
discretization-aware fit detects the excess dissipation
($\theta^2_h = -0.045$) and beats no-model uniformly. Notably, the upwind
errors \emph{decrease} with the CFL number: the energy injected by forward
Euler partially cancels the upwind over-dissipation---a fortuitous
cancellation between $r^h_{\mathrm{num}}$ and $r^\tau_{\mathrm{time}}$.
Supplying the temporal dissipation explicitly
therefore brings no benefit; the space--time closure is mildly
over-dissipative at the largest step ($0.078$ versus $0.055$ for
space-aware), the price of fitting cancelling parts independently.

\section{Conclusions and outlook}
\label{ag:sec:conclusion}
The forward-Euler time difference is a (one-sided) top-hat-filtered exact
time derivative, the temporal analogue of the spatial filter-swap.
Repeating the construction of
\cite{agdestein-jcp} in space \emph{and} time turns the fully discrete
forward-Euler finite-volume scheme into an exactly filtered conservation
law, whose residual flux
$r^{h, \tau} \mathrel{\mathop:}= \overline{f(u)}^{\tau} - f^h(\bar{u}^h)$ is a zero-error
a-posteriori target even at coarse time steps. The time integrator
contributes exactly one new term---a
flux time-quadrature error shrinking with the quadrature order---and no
temporal commutation term. A priori, this term
dominates the residual flux at practical CFL numbers; a posteriori,
a discretization-aware Smagorinsky model augmented with its leading
order---the Lax--Wendroff diffusion, dissipative by construction and linear in
$\tau$---keeps the error of a coarse forward-Euler simulation essentially
independent of the CFL number, where step-agnostic closures degrade. Much of
the time-quadrature error is non-dissipative and beyond the reach
of an eddy-viscosity form; machine-learned closures trained on the exact
target $r^{h, \tau}$, implicit and multi-stage integrators, and
three-dimensional Navier--Stokes are natural extensions.

\begin{acknowledgement}
This work is part of the project ``Discretize first, reduce next'' with
project number VI.Vidi.193.105 of the NWO Talent Programme Vidi, financed by
the Dutch Research Council (NWO). The author thanks Benjamin Sanderse and Roel
Verstappen for their supervision, and Josep Plana-Riu for discussions on the
dissipative properties of time integration.
\end{acknowledgement}


}

\backmatter
\appendix


\end{document}